\pdfoutput=1

\documentclass[final]{amsart}
\calclayout

\usepackage[utf8]{inputenc}
\usepackage{dsfont} 
\usepackage{comment}

\usepackage{hyperref}
\hypersetup{
	colorlinks = true,
	linkcolor = {blue},
	urlcolor = {black},
	citecolor = {black}
}

\usepackage{amssymb}
\usepackage{enumitem}
\usepackage[capitalize]{cleveref}
\usepackage{color}

\usepackage[backend=biber,
            style=alphabetic,
            bibencoding=utf8,
            language=auto,
            autolang=other,
            giveninits=true,
            doi=false,
            isbn=false,
            url=false,
            maxnames=10,
            sorting=nty]{biblatex}
\usepackage{asymptote}

\usepackage{mathtools}
\usepackage{faktor} 
\usepackage{tikz}
\usepackage{pgfplots}
\usepgfplotslibrary{patchplots,colormaps}
\pgfplotsset{compat = newest}
\usepackage{mathrsfs} 
\usetikzlibrary{babel} 
\usetikzlibrary{cd} 

\crefname{enumi}{}{} 

\newtheorem{theorem}{Theorem}[section]
\newtheorem*{theorem*}{Theorem}
\newtheorem{lemma}[theorem]{Lemma}
\newtheorem{proposition}[theorem]{Proposition}

\theoremstyle{remark}
\newtheorem{remark}{Remark}[section]
\newtheorem*{remark*}{Remark}

\theoremstyle{definition}
\newtheorem{definition}{Definition}[section]

\usepackage[backend=biber,style=alphabetic,bibencoding=utf8,language=auto,autolang=other,giveninits=true,doi=false,isbn=false,url=false,maxnames=10,sorting=nty]{biblatex}

\newcommand{\N}{\ensuremath{\mathbb N}} 
\newcommand{\R}{\ensuremath{\mathbb R}} 

\newcommand{\defeq}{\ensuremath{\coloneqq}}

\newcommand{\emptyparam}{\ensuremath{\,\cdot\,}}

\makeatletter 
\DeclarePairedDelimiter{\@tmpabs}{\lvert}{\rvert}
\newcommand{\@absstar}[1]{{\@tmpabs*{#1}}}
\newcommand{\@absnostar}[2][]{{\@tmpabs[#1]{#2}}}
\newcommand{\abs}{\@ifstar\@absstar\@absnostar}
\makeatother

\makeatletter 
\DeclarePairedDelimiter{\@tmpnorm}{\lVert}{\rVert}
\newcommand{\@normstar}[1]{{\@tmpnorm*{#1}}}
\newcommand{\@normnostar}[2][]{{\@tmpnorm[#1]{#2}}}
\newcommand{\norm}{\@ifstar\@normstar\@normnostar}
\makeatother

\newcommand\eps{\ensuremath{\varepsilon}}
\DeclareMathOperator{\Per}{Per}

\newcommand{\Haus}{\ensuremath{\mathscr H}} 

\usepackage{pythonhighlight} 

\DeclareMathOperator{\dist}{dist} 
\DeclareMathOperator{\diam}{diam}
\begin{document}

\title{Nonexistence of isoperimetric sets in spaces of positive curvature}

\author{Gioacchino Antonelli}
\address{Gioacchino Antonelli
\hfill\break Courant Institute Of Mathematical Sciences (NYU), 251 Mercer Street, 10012, New York, USA}
\email{ga2434@nyu.edu}

\author{Federico Glaudo}
\address{Federico Glaudo
\hfill\break School of Mathematics, Institute for Advanced Study, 1 Einstein Dr., Princeton NJ 05840, U.S.A.}
\email{fglaudo@ias.edu}

\begin{abstract}
For every $d\ge 3$, we construct a noncompact smooth $d$-dimensional Riemannian manifold with strictly positive sectional curvature without isoperimetric sets for any volume below $1$.
We construct a similar example also for the relative isoperimetric problem in (unbounded) convex sets in $\mathbb R^d$.
The examples we construct have nondegenerate asymptotic cone. 

The dimensional constraint $d\ge 3$ is sharp. 
Our examples exhibit nonexistence of isoperimetric sets only for small volumes; indeed in nonnegatively curved spaces with nondegenerate asymptotic cones isoperimetric sets with large volumes always exist.

This is the first instance of noncollapsed nonnegatively curved space without isoperimetric sets. 
\end{abstract}

\maketitle

\section{Introduction}
The aim of this paper is to prove the following two nonexistence results for isoperimetric sets (see \cref{def:isoperimetric}).

\begin{theorem}\label{thm:Main-boundary}
    For every $d\geq 3$, there exists a smooth complete $d$-dimensional Riemannian manifold $(M,g)$ with strictly positive sectional curvature everywhere and with nondegenerate asymptotic cone\footnote{Equivalently $\inf\limits_{p\in M,\, r>0}r^{-d}\mathrm{vol}_g(B_r(p))>0$.
    } for which the following hold.
    \begin{enumerate}
        \item If $v < 1$, there are no isoperimetric sets of volume $v$ in $M$;
        \item If $v>1$, there exists an isoperimetric set of volume $v$ in $M$.
    \end{enumerate}
\end{theorem}

\begin{theorem}\label{thm:Main}
    For every $d\geq 3$, there exists a closed strictly convex set $C\subseteq \mathbb R^d$ with nonempty interior, smooth boundary, and with nondegenerate asymptotic cone\footnote{Equivalently $\inf\limits_{p\in C,\, r>0}r^{-d}|B_r(p)|>0$.}
    for which the following hold.
    \begin{enumerate}
        \item If $v < 1$, there are no (relative) isoperimetric sets of volume $v$ in $C$;
        \item If $v>1$, there exists a (relative) isoperimetric set of volume $v$ in $C$.
    \end{enumerate}
\end{theorem}

These are the first instances of noncollapsed\footnote{A space $X$ is noncollapsed if $\inf_{x\in X} \abs{B_1(x)}>0$.} nonnegatively curved spaces without isoperimetric sets for some volume. \cref{thm:Main} answers a question raised in \cite[page 11]{LeonardiRitoreVernadakis2022}; \cref{thm:Main-boundary} and \cref{thm:Main} answer \cite[Question 1.5]{AP23} and \cite[Question 4.10]{PozzettaSurvey}. The hypothesis $d\geq 3$ in both \cref{thm:Main} and \cref{thm:Main-boundary} is sharp due to \cite{RitoreExistenceSurfaces01}, \cite[Theorem 1.4]{AP23}. Our examples exhibit nonexistence of isoperimetric sets only for \emph{small volumes}; indeed, in nonnegatively curved Alexandrov spaces (e.g., manifolds with nonnegative sectional curvature or convex sets in the Euclidean space) with nondegenerate asymptotic cone, isoperimetric sets with large volume always exist~\cite{LeonardiRitoreVernadakis2022, AntonelliBrueFogagnoloPozzetta2022, ConcavitySharpAPPS2}.

Let us remark that with our technique it should be possible to produce also examples of manifolds and convex sets satisfying \cref{thm:Main-boundary,thm:Main} with asymptotic cone that is \emph{not} nondegenerate.

\subsection{Overview of existence results for isoperimetric sets}
Our result complements the vast literature about the existence, under suitable assumptions, of isoperimetric sets in noncompact spaces (the existence in compact spaces is standard nowadays, see, e.g., the proof in \cite[Proposition 12.30]{Maggi2012}).
The usual assumptions are some form of nonnegativity of the (sectional, Ricci, or scalar) curvature and/or some rigidity of the structure of the space at infinity.

Under a (possibly negative) bound from below on the Ricci curvature, and a uniform lower bound on the volume of unit balls, one can show that, for a given volume, isoperimetric sets exist either in the space itself or the mass may split in a finite number of isoperimetric sets in its limits at infinity \cite{RitoreRosales, Nar14}. See also the generalizations to nonsmooth spaces in \cite{AntonelliNardulliPozzetta}. 

The \emph{weak existence} results mentioned in the previous paragraph are used in \cite{MondinoNardulli16} to deduce, e.g., the existence of isoperimetric sets for all volumes if the Ricci curvature is nonnegative and all the limits at infinity are Euclidean.

After the works in \cite{LeonardiRitoreVernadakis2022, AntonelliBrueFogagnoloPozzetta2022, ConcavitySharpAPPS2}, it is known that if the sectional curvature is nonnegative\footnote{Here it is sufficient to assume that the space is an Alexandrov space with nonnegative sectional curvature. Among the spaces satisfying this assumption we find: Riemannian manifolds with nonnegative sectional curvature, convex bodies in Euclidean spaces, and their boundaries.} and the asymptotic cone (cf., \cref{DefAsymptotic}) is nondegenerate, isoperimetric sets exist for all sufficiently large volumes. 
Analogous results, with slightly stronger assumptions on the asymptotic cones, are obtained also assuming only the nonnegativity of the Ricci curvature (in place of the nonnegativity of the sectional curvature). In the particular case of cones with nonnegative Ricci curvature, isoperimetric sets exist for all volumes and they are balls centered at one of the tips of the cone \cite{LionsPacella,MorganRitore02}.

The problem for $2$-dimensional spaces is completely settled in \cite{RitoreExistenceSurfaces01, AP23}: in nonnegatively curved surfaces existence holds for every volume.

The existence of isoperimetric sets for $3$-dimensional manifolds with nonnegative scalar curvature is known for asymptotically flat manifolds \cite[Proposition K.1]{Carlotto2016}, \cite{Shi2016}. The result holds also if one assumes that the manifold is asymptotic (in a quantitative sense) to a cone with nonnegative Ricci curvature \cite{ChodoshEichmairVolkmann17}, without any additional assumption on the curvature of the manifold. See also \cite{EichmairMetzger13,EichmairMetzger13b} for related results, and the recent survey \cite{BenattiFogagnolo2023}.

In the works \cite{LeonardiRitoreVernadakis2022,ConcavitySharpAPPS}, a number of techniques are developed to deal with the possible lack of isoperimetric sets for certain volumes to show the concavity of the isoperimetric profile for nonnegatively curved spaces (in the compact case this difficulty does not arise \cite[(iii) at page 483]{BavardPansu86},\cite{SternbergZumbrun1999,Kuwert,Bayle04}). The examples constructed in this paper tell us that these difficulties are unavoidable.

\subsection{Strategy of the construction and structure of the paper}
In \cref{sec:preliminaries} we recall the basic notions and some results about the isoperimetric problem in nonnegatively curved Alexandrov spaces.
The novelty of this section is \cref{lem:Threshold}: whenever the space has a nondegenerate asymptotic cone, if all limits at infinity are cones, then there is a threshold for the existence of isoperimetric sets (i.e., for $v<v_0$ there is no isoperimetric set with volume $v$ and for $v>v_0$ there is). To this end, we devise a robust proof of the concavity of the isoperimetric profile (to the power $\frac{d}{d-1}$) that might be of independent interest (see \cref{lem:decreasing} and \cref{rmk:robust-proof}).

In \cref{sec:construction}, that is the core of the paper, we construct the convex set $C$ which, after an approximation procedure performed in \cref{sec:approximation}, will satisfy the properties of \cref{thm:Main}. The manifold satisfying the properties mentioned in \cref{thm:Main-boundary} is obtained considering the boundary $\partial C$.

Let us quickly describe here the main features of our construction, we refer the reader to \cref{subsec:presentation-construction} for a more in-depth presentation and for a picture.

Consider a closed convex cone $\Sigma\subseteq\R^d$ contained in the half-space $\{x\in\R^d: x\cdot e_1\ge 0\}$. For a convex decreasing function $\varphi:\R\to(0,\infty)$ with $\varphi(+\infty)=0$, we define $C$ as 
\begin{equation*}
    C\defeq (\Sigma\times\R)\cap \{(x, t)\in\R^d\times\R:\, x\cdot e_1\ge \varphi(t)\}.
\end{equation*}
A rather delicate choice of $\varphi$ (so that $\varphi'$ is pointwise small compared to $\varphi$) produces a convex set $C$ without isoperimetric sets with small volume. The idea is that one limit at infinity of $C$ is the cone $\Sigma\times\R$, which has density strictly lower than any tangent cone at a point of $C$. This already implies that if there were isoperimetric sets for small volumes $v$ in $C$, then these sets would escape to infinity as $v\to 0$, see \cite[Theorem 6.9]{LeonardiRitoreVernadakis2022}.
The crux of the argument is to upgrade this observation to the nonexistence of isoperimetric sets for small volumes.

The proof is easier to follow for $C$ (i.e., for \cref{thm:Main}) than for $\partial C$ (i.e., for \cref{thm:Main-boundary}). To avoid repeating many arguments, we give the full proofs for $C$ and we emphasize only the differences in the proofs for $\partial C$.

Finally, in \cref{sec:proofs} we show how the results in the previous sections imply our two main theorems.

\subsection*{Acknowledgements}
The second author is supported by the National Science Foundation under Grant No. DMS--1926686.  

\section{Preliminaries}\label{sec:preliminaries}

\subsection{Alexandrov Spaces} 
We briefly introduce some key concepts about nonnegatively curved Alexandrov spaces. For the general theory we refer to \cite{BuragoBuragoIvanovBook,AKP}.

\subsubsection{General facts} 
A {\em nonnegatively curved Alexandrov space} is a complete, separable, geodesic space $(X,\mathrm{d})$ in which triangles are ``thicker'' than comparison triangles in $\mathbb R^2$, see, e.g., \cite[Definition 4.6.2]{BuragoBuragoIvanovBook}. Every Alexandrov space has either infinite Hausdorff dimension or integer Hausdorff dimension. We work only with finite-dimensional Alexandrov spaces, and we denote with $d\in\mathbb N$ the Hausdorff dimension of an Alexandrov space $X$. We denote with $\abs{\emptyparam}$ the $d$-dimensional Hausdorff measure $\Haus^d$ on a $d$-dimensional Alexandrov space. 

Notable examples of nonnegatively curved Alexandrov spaces are: complete Riemannian manifolds with nonnegative sectional curvature; closed convex sets with nonempty interior in the Euclidean space, endowed with the Euclidean distance; their boundaries endowed with their geodesic distance; metric cones on smooth compact Riemannian manifolds with sectional curvature $\geq 1$ (see \cite[Section 10.2]{BuragoBuragoIvanovBook}). In this paper we will be only concerned with closed convex sets with nonempty interior in the Euclidean space, and their boundaries, see \cite[Theorem 10.2.6]{BuragoBuragoIvanovBook} and \cite[Generalizations at p.359]{BuragoBuragoIvanovBook}. The theory of Alexandrov spaces provides a unified approach to address both settings. The reader not familiar with the theory of Alexandrov spaces can safely assume that whenever we state a result for a nonnegatively curved Alexandrov space, we will apply it only for convex sets or their boundaries.

Let us recall that if the boundary of a closed convex set with nonempty interior in $\mathbb R^d$ is smooth, then, when endowed with the metric tensor given by the pull-back of the Euclidean metric tensor, it is a smooth Riemannian manifold with nonnegative sectional curvature everywhere.

\subsubsection{Asymptotic cones and limits at infinity}

We are going to use the notion of pointed Gromov--Hausdorff (pGH) convergence (we refer the reader to \cite[Section 7 and Section 8.1]{BuragoBuragoIvanovBook}).
For the aims of this paper, since we are concerned only with Alexandrov spaces embedded in the Euclidean space, the notion of Hausdorff convergence on compact sets would be equivalent. The following statement formalizes this intuition.

\begin{lemma}\label{lem:boundary-GH-convergence}
    Given $d\ge 2$, let $(C_n)_{n\in\N}$ and $C$ be closed convex sets with nonempty interior in $\R^d$ such that $C_n$ converges to $C$ in the Hausdorff distance. Then $C_n$ converges to $C$ in the pointed Gromov-Hausdorff sense (for any sequence of points $C_n \ni p_n\to p\in C$).
    Moreover, the boundary $\partial C_n$, endowed with its geodesic distance, converges in the pointed Gromov-Hausdorff sense to the boundary $\partial C$ endowed with its geodesic distance.
\end{lemma}
\begin{proof}
    Upgrading the Hausdorff convergence to the pointed Gromov-Hausdorff convergence is simple and we leave it to the reader. The convergence at the level of the boundaries is proven in \cite[Theorem 1.2]{PETRUNIN97}.
\end{proof}

Let us define the notions of asymptotic cone and limits at infinity of a nonnegatively curved Alexandrov space $(X, \mathrm{d})$.

The \emph{asymptotic cone} $(X_\infty,\mathrm{d}_\infty)$ of $X$ is the limit (in the pointed Gromov-Hausdorff sense) 
\begin{equation}\label{DefAsymptotic}
(X,r_n\mathrm{d},o) \longrightarrow (X_\infty,\mathrm{d}_\infty,o_\infty),
\end{equation}
for some $r_n\to 0$, and for some arbitrary point $o\in X$.
The asymptotic cone does not depend on the choice of the rescaling factors $r_n$ or the point $o\in X$, it is itself a nonnegatively curved Alexandrov space, and it is a metric cone (see \cite[Definition 3.6.12]{BuragoBuragoIvanovBook} for the definition of metric cone). See, e.g., \cite[Exercises (a)--(e)]{BallmannGromovSchroeder85}.

For a $d$-dimensional nonnegatively curved Alexandrov space, we say that it has {\em nondegenerate asymptotic cone} if its asymptotic cone is $d$-dimensional.

We say that a metric space $(Y,\mathrm{d}_Y)$ is a {\em limit at infinity} of $X$ if there exists a sequence of diverging points $p_n\in X$, and $p_\infty\in Y$ such that 
\[
(X,\mathrm{d},p_n)\longrightarrow (Y,\mathrm{d}_Y,p_\infty).
\]

\subsubsection{Perimeter and isoperimetric profile}\label{sec:Perimeter}
As of today, a robust theory of functions of bounded variation and sets of finite perimeter has been developed in complete separable metric spaces endowed with a Radon measure \cite{Miranda03, AmbrosioDiMarino14}. Such a theory can be applied in the setting of $d$-dimensional nonnegatively curved Alexandrov spaces endowed with their Hausdorff measure $\Haus^d$, producing a well-behaved notion of {\em perimeter} of a set. We will not need the definitions here, but we refer the interested reader to \cite{Miranda03, AmbrosioDiMarino14} for the general theory.

Given a nonnegatively curved Alexandrov space $X$ and a Borel set $E\subset X$, we denote by $\Per_X(E)$ the perimeter of $E$ in $X$, defined as in \cite[Definition 4.1]{Miranda03}. This notion of perimeter agrees with the usual one in the following cases: 
\begin{enumerate}
    \item Assume that $X$ is a closed convex set with nonempty interior in $\mathbb R^d$, for $d\ge 2$. For every Borel set $E\subset X$ with $\Per_X(E)<\infty$, we have that $E$ is a set of finite perimeter and, denoting with $\partial^* E$ its reduced boundary, $\Per_X(E)=\Haus^{d-1}(\partial^* E\cap\mathrm{int}(X))$. Namely $\Per_X(E)$ is the {\em relative} perimeter in the interior $\mathrm{int}(X)$.
    \item Assume that $X$ is a smooth $d$-dimensional Riemannian manifold, for $d\ge 2$. For every Borel set $E\subset X$ with $\Per_X(E)<\infty$, we have that $E$ is a set of finite perimeter and, denoting with $\partial^*E$ its reduced boundary, $\Per_X(E)=\Haus^{d-1}(\partial^* E)$.
\end{enumerate}

\begin{definition}\label{def:isoperimetric}
For a nonnegatively curved Alexandrov space $X$, we denote by $I_X$ its {\em isoperimetric profile}, i.e., the function
\[
(0,+\infty)\ni v\mapsto I_X(v):=\inf\big\{\Per_X(E):\,\abs{E}=v\big\}.
\]
We say that a Borel set $E\subseteq X$ is an {\em isoperimetric set} if $\Per_X(E)=I_X(|E|)$. We will always implicitly assume that isoperimetric sets are open~\cite[Theorem 1.4]{AntonelliPasqualettoPozzetta}.
\end{definition}

\subsection{Known results on the isoperimetric profile and isoperimetric sets}
Let us recall some standard facts in the theory of the isoperimetric problem for nonnegatively curved spaces.\\
First of all, by scaling-invariance, if $X$ is a $d$-dimensional metric cone then the isoperimetric profile satisfies $I_X(v) = I_X(1)v^{\frac{d-1}d}$ (see \cite[Theorem 1.1]{LionsPacella}, \cite[Proposition 3.1]{RitoreRosales}, or \cite[Theorem 3.6]{MorganRitore02} for the sharp isoperimetric inequality in convex cones in the Euclidean space, in Euclidean cones, or in nonnegatively Ricci curved cones). From now on, we will use the latter information without mentioning it.

The following lemma compares the isoperimetric profile of a nonnegatively curved Alexandrov space with the one of its asymptotic cone and limits at infinity.

\begin{proposition}\label{prop:profile-properties}
    Given $d\ge 2$, let $X$ be a $d$-dimensional nonnegatively curved Alexandrov space with nondegenerate asymptotic cone $X_\infty$.
    The function $v\mapsto I_X(v)^{\frac{d}{d-1}}$ is concave. Moreover, for all $v>0$, we have
    \begin{equation*}
        0 < I_{X_{\infty}}(v) \le I_X(v) \le I_{X'}(v),
    \end{equation*}
    for any $X'$ limit at infinity of $X$. If $I_{X_\infty}(v)=I_X(v)$ for some $v>0$ then $X$ is isometric to $X_\infty$.
\end{proposition}
\begin{proof}
    From \cite[Theorem 1.1]{ConcavitySharpAPPS} we get that $v\mapsto I_{X}^{\frac{d}{d-1}}(v)$ is a concave function.\\
    The inequality $I_X\leq I_{X'}$, for any limit at infinity $X'$, follows from the more general \cite[Equation (2.17) in Proposition 2.19]{AntonelliNardulliPozzetta}, see also \cite[Proposition 4.3]{LeonardiRitoreVernadakis2022} in the setting of convex bodies.\\
    The inequality $I_{X_\infty}\leq I_X$ follows from the more general \cite[Theorem 1.1]{BaloghKristaly}, see also \cite[Theorem 6.3]{LeonardiRitoreVernadakis2022} in the setting of convex bodies.
    The last part of the statement follows from the rigidity statement in \cite[Item (3) of Theorem 1.2]{ConcavitySharpAPPS2}. 
\end{proof}

The following lemma is a direct consequence of \cite[Lemma 4.20]{ConcavitySharpAPPS}, see also \cite[Theorem 5.8]{LeonardiRitoreVernadakis2022} in the case of convex bodies.
\begin{lemma}\label{lem:LimitiVannoAInfinito}
    Given $d\ge 2$, let $X$ be a $d$-dimensional nonnegatively curved Alexandrov space with nondegenerate asymptotic cone.
    For all $v>0$ the following holds: if there is no isoperimetric set with volume $v$ in $X$, then there exists an isoperimetric set of volume $v$ in a limit at infinity $X'$ of $X$, and moreover $I_{X'}(v)=I_X(v)$.
\end{lemma}

Let us recall the following standard result on the diameter of isoperimetric sets.
See \cite[Lemma 5.5]{LeonardiRitoreVernadakis2022} for the version on convex sets, and \cite[Proposition 4.23]{ConcavitySharpAPPS2} for the general version on nonnegatively curved Alexandrov spaces.

\begin{lemma}\label{lem:Boundedness}
    Given $d\geq 1$, let $X$ be a $d$-dimensional nonnegatively curved Alexandrov space with nondegenerate asymptotic cone, and assume that $\abs{B_1^{X_\infty}} \ge \vartheta_0>0$, where $B_1^{X^\infty}\subseteq X_\infty$ denotes the unit ball centered at a tip of the asymptotic cone $X_\infty$. 
    There exists a constant $A_2:=A_2(d, \vartheta_0)>0$ such that every isoperimetric set $E$ in $X$
    satisfies the following inequality
    \[
    \mathrm{diam}(E)\leq A_2|E|^{\frac{1}{d}},
    \]
    where $\mathrm{diam}(E)$ denotes the diameter of $E$.
\end{lemma}

\subsection{The existence of a threshold}
The goal of this section is to show that, under a mild assumption on the limits at infinity, there is a volumetric threshold for the existence of isoperimetric sets, see \cref{lem:Threshold}.

Our main tool is the following \emph{robust} proof of the concavity of $I_X^{\frac{d}{d-1}}$ under the assumption of existence of an isoperimetric set. We postpone to \cref{rmk:robust-proof} a comparison between \cref{lem:decreasing} and the standard proof of the concavity of $I_X^{\frac{d}{d-1}}$.
\begin{lemma}\label{lem:decreasing}
    Given $d\geq 2$, let $X$ be a nonnegatively curved $d$-dimensional Alexandrov space. Let $E\subseteq X$ be an isoperimetric set. For $r\in\R$, define $E_r$ as\footnote{We denote with $\mathrm{d}_F:X\to[0,\infty)$ the distance from the set $F$.}
    \begin{equation*}
        E_r \defeq
        \begin{cases}
            \{x\in X: \mathrm{d}_{\overline E}(x)\leq r\}
            &\quad\text{if $r \ge 0$,} \\
            \{x\in X: \mathrm{d}_{X\setminus E}(x)\ge -r\} 
            &\quad\text{if $r < 0$}.
        \end{cases}
    \end{equation*}
    There exists $H\geq 0$, which agrees with the (constant) mean curvature of $\partial E$ whenever $X$ is sufficiently smooth and which agrees with $I_X'(\abs{E})$ if $I_X$ is differentiable at $\abs{E}$,
    such that the function
    \begin{equation}\label{Bellaaa}
        \R\ni r\mapsto \Per_X(E_r)^{\frac{d}{d-1}}-\frac{d}{d-1}H\Per_X(E)^{\frac{1}{d-1}}\abs{E_r},
    \end{equation}
    achieves its maximum at $r=0$ on $(-\infty,\infty)$. Moreover, for $r\ge 0$, we have
    \begin{equation}\label{eq:enlargement-ratio}
        \frac{\Per_X(E_r)^{\frac{d}{d-1}}}{\abs{E_r}} \le 
        \frac{\Per_X(E)^{\frac{d}{d-1}}}{\abs{E}}.
    \end{equation}
\end{lemma}
\begin{proof}
    First of all, the existence of an isoperimetric set $E$ implies that $\inf_{x\in X} \abs{B_1(x)}>0$ (see~\cite[Proposition 2.18]{AntonelliBrueFogagnoloPozzetta2022}); this assumption is necessary in some of the results we will cite.
    We will use repeatedly, without saying it explicitly, that the map $r\mapsto\Per_X(E_r)$ is lower semicontinuous everywhere. In addition, by the co-area formula \cite[Proposition 4.2]{Miranda03}, we also have that $\frac{\mathrm d}{\mathrm dr}\abs{E_r}=\Per_X(E_r)$ for almost every $r\in \R$ and in the distributional sense, and $\Per_X(E_r)\in L^1_{\mathrm{loc}}(\R)$.
    
    Let $H \geq 0$ be a mean curvature barrier for $E$, as defined in \cite[Definition 3.6]{ConcavitySharpAPPS}. The existence of such an $H$ is guaranteed by \cite[Theorem 3.3]{ConcavitySharpAPPS}, and the fact that $H\geq 0$ is guaranteed by \cite[Proposition 3.1]{ConcavitySharpAPPS2}, since $X$ has nonnegative curvature. The fact that $H$ agrees with the usual mean curvature in the smooth case is observed in \cite[Remark 3.8]{ConcavitySharpAPPS}. Finally, the fact that $H=I'_X(|E|)$, if $I_X$ is differentiable at $|E|$, comes from \cite[Corollary 4.21]{ConcavitySharpAPPS}.
    
    From \cite[Equation (3.37)]{ConcavitySharpAPPS} we get 
    \begin{equation}\label{eqn:Step2}
        \Per_X(E_r)\leq \Per_X(E)\left(1+\frac{H}{d-1}r\right)^{d-1} \quad\text{for all $r> -\frac{d-1}{H}$},
    \end{equation}
    and $\Per_X(E_r)=0$ for every $r\leq -\frac{d-1}{H}$. When $H=0$ we understand that $-\frac{d-1}{H}=-\infty$. Thus, since $r\mapsto \Per_X(E_r)$ is lower semicontinuous, we get that $r\mapsto \Per_X(E_r)$ is continuous at $r=0$. We will use this information repeatedly without saying it. Notice that, as a byproduct of \cref{eqn:Step2}, we also get that $\Per_X(E_r)$ is locally bounded.

    Following the proof of \cite[Equation (3.34)]{ConcavitySharpAPPS}, but integrating on the slab $E_R\setminus E_r$ with $R\geq r\geq 0$ instead of the slab $E_r\setminus E$, we obtain 
    \begin{equation}\label{eqn:KeyEstimate}
    \begin{aligned}
        &\Per_X(E_R)\leq \Per_X(E_r) +\int_r^R\frac{H}{1+\frac{H}{d-1}\varrho}\Per_X(E_\varrho)\mathrm{d}\varrho
        &\quad \text{for almost all $0 \le r \le R$.}
    \end{aligned}
    \end{equation}
    From \cref{eqn:KeyEstimate} we can infer
    \begin{equation}\label{eqn:Step1}
    \begin{aligned}
        \frac{\mathrm{d}}{\mathrm{d}r}\Per_X(E_r)\leq \Per_X(E_r)\frac{H}{1+\frac{H}{d-1}r} 
        &\quad\text{distributionally on $(0,\infty)$.}
    \end{aligned}
    \end{equation}
    
    Let us now show that the function in the statement is distributionally (weakly) decreasing for $r\in [0,+\infty)$.
    For $r\in(0,\infty)$, it holds in the distributional sense that
    \begin{equation}
        \begin{split}
            \frac{\mathrm{d}}{\mathrm{d}r}&\left(\Per_X(E_r)^{\frac{d}{d-1}}-\frac{d}{d-1}H\Per_X(E)^{\frac{1}{d-1}}|E_r|\right)\\
            &=\frac{d}{d-1}\Per_X(E_r)\left(\Per_X(E_r)^{\frac{1}{d-1}}\frac{\frac{\mathrm{d}}{\mathrm{d}r}\Per_X(E_r)}{\Per_X(E_r)}-H\Per_X(E)^{\frac{1}{d-1}}\right)\\
            &\leq \frac{dH}{d-1}\Per_X(E_r)\left(\Per_X(E_r)^{\frac{1}{d-1}}\frac{1}{1+\frac{H}{d-1}r}-\Per_X(E)^{\frac{1}{d-1}}\right)\leq 0,
        \end{split}
    \end{equation}
    where in the first inequality we used \cref{eqn:Step1}, and in the second inequality we used \cref{eqn:Step2}. 
    We deduce that the function in the statement has a maximum at $r=0$ when restricted to $r\in [0, \infty)$.

    Let us consider the case $r\in (-(d-1)/H,0)$. This is enough since the function in the statement is constantly zero on $(-\infty,-(d-1)/H]$. We have, by using \cref{eqn:Step2}, 
    \begin{equation}\label{NartroPasso}
        \begin{split}
            \frac{d}{d-1}&H\left(|E|-|E_r|\right) = \frac{d}{d-1}H\int_r^0\Per_X(E_\varrho)\mathrm{d}\varrho\\
            &\leq \frac{d}{d-1}H\Per_X(E)\int_r^0 \left(1+\frac{H}{d-1}\varrho\right)^{d-1}\mathrm{d}\varrho \\
            &= \frac{d}{d-1}H\Per_X(E)\left[\frac{d-1}{Hd}\left(1-\left(1+\frac{H}{d-1}r\right)^d\right)\right]\\
            &=\Per_X(E)\left(1-\left(1+\frac{H}{d-1}r\right)^d\right).
        \end{split}
    \end{equation}
    Thus, by subsequently using \cref{eqn:Step2} and \cref{NartroPasso} we get 
    \begin{equation}
        \begin{split}
            \Per_X(E_r)^{\frac{d}{d-1}}&-\frac{d}{d-1}H\Per_X(E)^{\frac{1}{d-1}}|E_r|\\
            &\leq\Per_X(E)^{\frac{1}{d-1}}\left[\Per_X(E)\left(1+\frac{H}{d-1}r\right)^d-\frac{d}{d-1}H\abs{E_r}\right]\\
            &\leq \Per_X(E)^{\frac{1}{d-1}}\left(\Per_X(E)-\frac{d}{d-1}H|E|\right), 
        \end{split}
    \end{equation}
    thus the function in the statement has a maximum at $r=0$ when restricted to $r\in (-\infty,0]$.

    It remains to prove \cref{eq:enlargement-ratio}. Using \cref{Bellaaa}, we have
    \begin{equation*}
        \frac{\Per_X(E_r)^{\frac{d}{d-1}}-\Per_X(E)^{\frac{d}{d-1}}}
        {\abs{E_r}-\abs{E}}
        \le \frac{d}{d-1}H\Per_X(E)^{\frac1{d-1}}
        \le
        \frac{\Per_X(E)^{\frac{d}{d-1}}}
        {\abs{E}}
    \end{equation*}
    and, for $r>0$, the inequality \cref{eq:enlargement-ratio} follows.
\end{proof}

\begin{remark}
    The statement and proof of \cref{lem:decreasing} work in the more general setting of noncompact metric measure spaces with synthetic nonnegative Ricci curvature. For the ease of readability and for consistency with the material of this paper we stated it for nonnegatively curved Alexandrov spaces.
\end{remark}

\begin{remark}\label{rmk:robust-proof}
    Let $X$ be a $d$-dimensional nonnegatively curved Alexandrov space.
    Let us see how \cref{lem:decreasing} implies the concavity of $I_X^{\frac{d}{d-1}}$ and why it is a more robust method compared to the standard proof (see, e.g., \cite[Theorem 4.4]{ConcavitySharpAPPS}, or \cite[Theorem 2.1]{Bayle04} in the smooth setting).
    The proof starts by assuming, for each $v_0>0$, the existence of an isoperimetric set $E$ with $\abs{E}=v_0$ (possibly looking in a limit at infinity, and possibly considering simultaneously multiple limits at infinity, thus using~\cite[Theorem 1.1]{AntonelliNardulliPozzetta}).

    Adopting the same assumptions and notation of \cref{lem:decreasing}, consider the function $(0,\infty)\ni \abs{E_r} \mapsto J(\abs{E_r})\defeq \Per_X(E_r)^{\frac d{d-1}}$.
    By definition of isoperimetric profile, we have $I_X(v)^{\frac d{d-1}}\le J(v)$ for all $v>0$ and $I_X(v_0)^{\frac d{d-1}}=J(v_0)$. Moreover, the statement of \cref{lem:decreasing} is equivalent to $J(v) \le J(v_0) + (v-v_0)\gamma$, where $\gamma\defeq \frac{d}{d-1} H\Per_X(E)^{\frac 1{d-1}}$. Therefore, we obtain
    \begin{equation*}
        I_X(v)^{\frac{d}{d-1}} \le J(v) \le J(v_0)+(v-v_0)\gamma = I_X(v_0)^{\frac{d}{d-1}} + (v-v_0)\gamma,
    \end{equation*}
    that is, we have found a line that touches $I_X(v)^{\frac{d}{d-1}}$ at $v=v_0$ and stays above (or equal) for all other values. The concavity of $I_X^{\frac{d}{d-1}}$ follows.

    The standard proof goes along the same scheme but replaces $J(v)\le J(v_0)+(v-v_0)\gamma$ with the weaker (and \emph{differential} in nature) $J(v)\le J(v_0)+(v-v_0)\gamma + o(\abs{v-v_0})$ which is sufficient for the concavity but not for the rigidity argument employed in the proof of \cref{lem:Threshold}.
\end{remark}

\begin{lemma}\label{lem:Threshold}
    Given $d\geq 1$, let $X$ be a $d$-dimensional nonnegatively curved Alexandrov space with nondegenerate asymptotic cone. Assume that all its limits at infinity are cones.

    Then there exists $0 \le v_0<+\infty$ such that the following hold.
    \begin{enumerate}
        \item For $0<v< v_0$, there are no isoperimetric sets of volume $v$ in $X$;
        \item For $v_0 < v$, there exists an isoperimetric set with volume $v$ in $X$.
    \end{enumerate}
\end{lemma}
\begin{proof}
    If isoperimetric sets exist for every positive volume, $v_0=0$ works. From now on we assume that this is not the case.

    Let
    \begin{equation*}
        \vartheta\defeq\inf \{I_{X'}(1): \, X' \text{ is a limit at infinity of $X$}\}
    \end{equation*}
    be the infimum of the perimeter of sets with unit volume in a limit at infinity of $X$. Denote with $X_\infty$ the asymptotic cone of $X$. For all $v>0$, as a consequence of \cref{prop:profile-properties}, we have
    \begin{equation}\label{eq:profiles-ordered}
        0 < I_{X_\infty}(v) \le I_X(v) \le \vartheta v^{\frac{d-1}{d}} \le I_{X'}(v),
    \end{equation}
    for any $X'$ limit at infinity of $X$ (we are using that all the limits at infinity are cones).
    In particular $\vartheta > 0$.

    Observe that $I_X(v)/v^{\frac{d-1}{d}}$ is a decreasing continuous quantity that is $\leq \vartheta$ everywhere; let $v_0\in [0,+\infty]$ be the maximum value such that $I_X(v_0)=\vartheta v_0^{\frac{d-1}{d}}$. We show that the statement holds for this choice of $v_0$. Notice that a posteriori $v_0<+\infty$ because we always have existence of isoperimetric sets for large volumes, see \cite[Item (1) of Theorem 1.2]{ConcavitySharpAPPS2}.

    Assume that there is no isoperimetric set of volume $\tilde v>0$ in $X$.
    Then there exists an isoperimetric set with volume $\tilde v$ and perimeter $I_X(\tilde v)$ in one of the limits at infinity of $X$, see \cref{lem:LimitiVannoAInfinito}. Therefore, taking \cref{eq:profiles-ordered} into account, we get $I_X(\tilde v) = \vartheta \tilde v^{\frac{d-1}{d}}$. In particular $\tilde v\le v_0$.

    Take $0<\tilde v<v_0$ and assume by contradiction that there is an isoperimetric set $E\subseteq X$ of volume $\tilde v$ with perimeter equal to $I_X(\tilde v)=\vartheta \tilde v^{\frac{d-1}{d}}$. Since $I_X(v)=\vartheta v^{\frac{d-1}{d}}$ in a neighborhood of $\tilde v$, the constant $H$ appearing in the statement of \cref{lem:decreasing} coincides with $H=I'_X(\tilde v) = \frac{d-1}{d}\vartheta \tilde v^{-\frac1d}$ and so $\frac{d}{d-1}H\Per_X(E)^{\frac1{d-1}} = \vartheta^{\frac{d}{d-1}}$. For any $r\in\R$, let $E_r$ be the set defined in \cref{lem:decreasing}; we have
    \begin{equation*}
        \Per_X(E_r)^{\frac{d}{d-1}}
        -
        \vartheta^{\frac{d}{d-1}}\abs{E_r}
        \le 
        \Per_X(E)^{\frac{d}{d-1}}
        -
        \vartheta^{\frac{d}{d-1}}\abs{E}
        = 0.
    \end{equation*}
    Hence, for any $v>0$, since one can find $r\in\R$ so that $\abs{E_r}=v$, there is a subset $E'\subseteq X$ with $\abs{E'}=v$ so that
    \begin{equation*}
        \Per_X(E') \le \vartheta v^{\frac{d-1}{d}}.
    \end{equation*}
    This is a contradiction because we are assuming (see the beginning of the proof) that there is at least one volume $v\in\R$ so that there is not an isoperimetric set of volume $v$ and we have shown above that in such case it must hold $I_X(v)=\vartheta v^{\frac{d-1}{d}}$.
\end{proof}
\begin{remark}
    The authors believe that the assumption that the limits at infinity are cones is necessary in the statement of the previous lemma. In particular, it might be possible to construct a counterexample to the statement with the same method of this paper (choosing as $\Sigma$ in \cref{eq:sigma-properties} a convex set that is a compact perturbation of a cone).

    For $v=v_0$, under the same assumptions of \cref{lem:Threshold}, we expect an isoperimetric set to exist in certain cases and not to exist in others.
\end{remark}

\section{The construction}\label{sec:construction}

\subsection{Presentation of the construction}\label{subsec:presentation-construction}
Let $d\geq 2$ be a natural number.
Let $\Sigma\subseteq\R^d$ be a subset such that
\begin{equation}\label{eq:sigma-properties}
    \begin{aligned}
        &\Sigma \text{ is a closed convex cone with nonempty interior,}\\
        &\Sigma\setminus\{0_{\R^d}\}\subseteq \{x\in\R^d:\, x\cdot e_1 > 0\}, \\
        & \text{The boundary $\partial\Sigma$ is smooth away from the origin.}
    \end{aligned}
\end{equation}

Let $\varphi:\R\to(0,\infty)$ be a function such that
\begin{equation}\label{eq:phi-properties}
    \begin{aligned}
        &\varphi \text{ is convex and decreasing,}\\
        &\lim_{t\to+\infty}\varphi(t) = 0,\\
        &\varphi \text{ is affine on $(-\infty, 0]$ (i.e., $\varphi''=0$ on $(-\infty, 0]$).}
    \end{aligned}
\end{equation}
Let $C\subseteq\R^{d+1}$ be the closed convex set with nonempty interior (here $e_1$ denotes the first element of the canonical basis of $\R^d$)
\begin{equation}\label{eq:def-C}
    C \defeq (\Sigma\times\R) \cap \{(x, t)\in \R^d\times \R:\, x \cdot e_1 \geq \varphi(t)\}.
\end{equation}

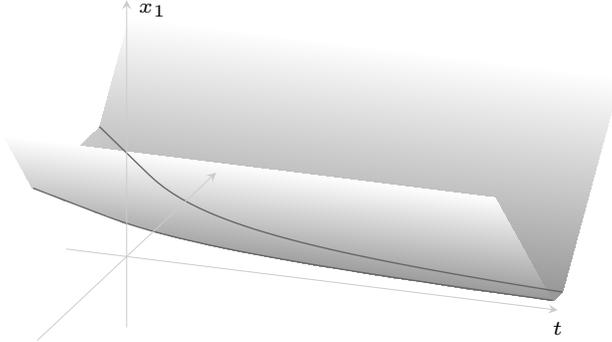
\begin{figure}[htbp]\label{fig}

 
\begin{tikzpicture}[scale=1.3,
declare function={
phi(\x)= (\x < 0) * (1/1.412-0.5*\x)   +
         (\x >= 0) * (1/(\x+1.412))
;
}]
 
\begin{axis} [
    axis lines = center, 
    axis on top, 
    axis line style = {very thin, black!20},
    view = {20}{20},
    unit vector ratio = 1 1 1,
    xtick = {0},
    ytick = {0},
    ztick = {0},
    xlabel = {{\tiny $t$}},
    zlabel = {{\tiny $x_1$}},
    x label style={anchor=north},
]
 
\addplot3 [
    domain=-0.7:5,
    domain y = -2:2,
    samples = 20,
    samples y = 100,
    surf,
    shader = interp,
    opacity = 0.4,
    colormap/gray] {max(phi(x), abs(y))};
    
\addplot3+[
    black!60,
    no markers,
    samples=51, 
    samples y=0,
    domain=-0.7:5,
    variable=\x]
    ({\x},{phi(\x)},{phi(\x)});

\addplot3+[
    black!60,
    no markers,
    samples=51, 
    samples y=0,
    domain=-0.7:5,
    variable=\x]
    ({\x},{-phi(\x)},{0.01+phi(\x)});

\end{axis}
 
\end{tikzpicture}
 
    \caption{Representation of the boundary $\partial C$. Notice that any slice of $C$ at a fixed $t=t_0$ coincides with $\Sigma_\lambda$ for some $\lambda>0$.}
\end{figure}

\subsection{Properties of \texorpdfstring{$C$}{C}}\label{PropertiesC}
The asymptotic cone of $C$ is
\begin{equation*}
    (\Sigma\times\R)\cap \{(x, t):\, x\cdot e_1 \geq \varphi'(0)t\}.
\end{equation*}
Notice that the limit at infinity of $C$ along the sequence of points $(\varphi(t)e_1, t)$ for $t\to\infty$ is $\Sigma\times\R$. As a consequence it holds that $I_{\Sigma\times\R}\ge I_C$ (recall \cref{prop:profile-properties}). Finally notice that all the limits at infinity of $C$ are cones.

For $h\in\R$, let
\begin{equation*}
    \Sigma_{h} \defeq \Sigma\cap \{x\in\R^d: x\cdot e_1 \geq h\} = h\Sigma_1.
\end{equation*}
Notice that $\Sigma_{\varphi(t_0)}$ is the section of $C$ at $t=t_0$.

Let $A_1=A_1(\Sigma)\defeq \sup_{x\in\Sigma}\frac{\abs{x}}{x\cdot e_1}<\infty$.
Observe that in the region $\{(x, t): \abs{x}>A_1\varphi(0)\}$, the set $C$ coincides with the cone
\begin{equation}\label{eq:outside-affine-cone}
    C'\defeq (\Sigma\times\R)\cap \{(x, t): x\cdot e_1 \geq \varphi(0) + t\varphi'(0)\} ,
\end{equation}
the tip of the cone being $\left(0_{\R^d}, -\frac{\varphi(0)}{\varphi'(0)}\right)$. One can show this by looking at the sections of $C$ and $C'$. For $t_0<0$, the sections of $C$ and $C'$ at $t=t_0$ coincide.
To handle the sections for $t_0\ge 0$, notice that if $(x, t)\in(\Sigma\times\R)\cap \{(x, t): \abs{x}\geq A_1\varphi(0)\}$ then $x\cdot e_1 \ge \varphi(0)\ge \varphi(t_0)$. Therefore, the sections of $C$ and $C'$ at $t\ge0$ coincide with $\Sigma$ when we consider only the region $\{x: \abs{x}>A_1\varphi(0)\}$.

Let us remark that $\Sigma_1\times\R$ is not isometric to $\Sigma\times\R$. One can show this observing that the density of $\Sigma_1$ at any of its points is strictly larger than the density of $\Sigma$ at its tip; thus the density of $\Sigma_1\times\R$ at any point is strictly larger than the density of $\Sigma\times\R$ at its tip.

\subsection{Properties of the boundary \texorpdfstring{$\partial C$}{of C}}
Let us now shift our attention to the boundary of $C$. Since it is the boundary of a convex set, it is in particular an Alexandrov space with nonnegative sectional curvature when endowed with its geodesic distance. 
The asymptotic cone of $\partial C$ is the boundary of the asymptotic cone of $C$ (see \cref{lem:boundary-GH-convergence}). 
Whenever we consider the boundary of a convex set, we implicitly assume that it is endowed with the geodesic distance.

Notice that the limit at infinity of $\partial C$ along the sequence of points $(\varphi(t)e_1, t)$ for $t\to\infty$ is $\partial\Sigma\times\R$ (see again \cref{lem:boundary-GH-convergence}).
As a consequence it holds that $I_{\partial\Sigma\times\R}\ge I_{\partial C}$ (recall \cref{prop:profile-properties}).
The section of $\partial C$ at $t=t_0$ is $\partial\Sigma_{\varphi(t_0)}$.
Analogously to what we observed above for $C$, also $\partial C$ coincides with a cone (that is $\partial C'$) in the region $\{(x, t): \abs{x}\geq A_1\varphi(0)\}$.

For $d\ge 3$, we have that $\partial\Sigma_1\times\R$ is not isometric to $\partial\Sigma\times\R$ (notice that for $d=2$, $\partial\Sigma_1$ and $\partial\Sigma$ are both isometric to $\R$).  This can be shown by observing that the tangent cone at any point of $\partial\Sigma_1$ splits at least one line\footnote{Here we are using that $\partial\Sigma$ is smooth away from the origin. We need this assumption only here. The fact that $\partial\Sigma_1\times\R$ is not isometric to $\partial\Sigma\times\R$ should hold also without the smoothness assumption, but the proof becomes more involved.}, while the tangent cone at the tip of $\partial\Sigma$ does not split any line.

\subsection{Diameter bound for isoperimetric sets in \texorpdfstring{$C$ and $\partial C$}{C and its boundary}}
The asymptotic cone of $C$ contains $\Sigma\times[0,\infty)$. The asymptotic cone of $\partial C$ contains $\partial\Sigma\times[0,\infty)$. These two remarks imply a lower bound on the volumes of the unit balls of the asymptotic cones of $C$ and $\partial C$ that is independent of the choice of $\varphi$. Therefore \cref{lem:Boundedness} tells us that there exist two constants $A_2=A_2(\Sigma), A_2^\partial=A_2^\partial(\Sigma)$ such that any isoperimetric set $E\subseteq C$ satisfies
\begin{equation}\label{eq:diam-C}
    \diam(E) \le A_2 \abs{E}^{\frac{1}{d+1}},
\end{equation}
and any isoperimetric set $E\subseteq\partial C$ satisfies
\begin{equation}\label{eq:diam-C-boundary}
    \diam(E) \le A_2^\partial\Haus^d(E)^{\frac{1}{d}}.
\end{equation}

\subsection{Controlling the isoperimetric profile of \texorpdfstring{$C$ and $\partial C$}{C and its boundary}}
The following two lemmas (one for $C$ and one for $\partial C$) provide a lower bound on the perimeter of a subset that is localized in space.

\begin{lemma}\label{lem:crucial}
    Given $d\ge 2$, let $\Sigma\subseteq\R^d$ be a subset satisfying \cref{eq:sigma-properties}, 
    let $\varphi:\R\to(0,\infty)$ be a function satisfying \cref{eq:phi-properties}, and let $C\subseteq\R^{d+1}$ be the set defined in \cref{eq:def-C}.

    Given $R>0$ and $\alpha>0$, there exists $A_3=A_3(d, R, \alpha)\in (0, 1)$ such that the following statement holds.

    Given $t_0\in\R$, consider a Borel set $E\subseteq C\cap (B_R(0_{\R^d})\times[t_0, t_0+R])$ with $\abs{E}=1$.
    If $\varphi(t_0)\ge \alpha$ and $\abs{\varphi'(t_0)}\le A_3$, then
    \begin{equation*}
        \Per_C(E) > I_{\Sigma\times\R}(1).
    \end{equation*}
\end{lemma}
\begin{proof}
    By translation invariance in the $t$ variable\footnote{To be precise the assumptions on $\varphi$ (see \cref{eq:phi-properties}) are not entirely translation invariant because we assume that $\varphi$ is affine on $(-\infty, 0]$. This assumption will not play any role in this proof.} we can assume $t_0=0$.
    Furthermore, since $C\cap (B_R(0_{\R^d})\times[0, R])$ is nonempty, it must hold $\varphi(R)\le R$.

    Let us argue by contradiction. Assume that there is a sequence of functions $\varphi_n:\R\to[0,\infty)$, which generate $C_n$ as in \cref{eq:def-C}, and sets $E_n\subseteq C_n\cap (B_R(0_{\R^d})\times[0, R])$ with $\abs{E_n}=1$ so that
    \begin{align*}
        &\alpha\le \varphi_n(0), \varphi_n(R)\le R ,\\
        &\varphi_n'(0)\to 0 \text{ as $n\to\infty$},\\
        &\Per_C(E_n) \le I_{\Sigma\times\R}(1) .
    \end{align*}
    Observe that the derivatives are going to $0$ uniformly in $[0,R]$ (because of the convexity of $\varphi_n$) and therefore we may assume that, up to extracting a subsequence, $\varphi_n\to\lambda\in[\alpha, R]$ uniformly in $[0, R]$.
    Hence, the spaces $C_n\cap (B_R(0_{\R^d})\times [0, R])$ converge in the Gromov-Hausdorff topology to $(\Sigma_\lambda\times\R)\cap (B_R(0_{\R^d})\times [0, R])$. Moreover the sets $E_n$ converge to a set $E\subseteq \Sigma_\lambda\times\R$ with $\abs{E}=1$, and thus by lower semicontinuity of the perimeter we have $\Per_{\Sigma_\lambda\times\R}(E)\le I_{\Sigma\times\R}(1)$ (see \cite[Proposition 3.6]{AmbrosioBrueSemola19}). This is a contradiction because, thanks to \cref{prop:profile-properties} (we are using that $\Sigma_\lambda\times\R$ is not isometric to $\Sigma\times\R$), we have $\Per_{\Sigma_\lambda\times\R}(E) \ge I_{\Sigma_\lambda\times\R}(1) > I_{\Sigma\times\R}(1)$.  
\end{proof}

Now we state the equivalent lemma on the boundary $\partial C$.
\begin{lemma}\label{lem:crucial-boundary}
    Given $d\ge 3$, let $\Sigma\subseteq\R^d$ be a subset satisfying \cref{eq:sigma-properties}, 
    let $\varphi:\R\to(0,\infty)$ be a function satisfying \cref{eq:phi-properties}, and let $C\subseteq\R^{d+1}$ be the set defined in \cref{eq:def-C}.

    Given $R>0$, $\alpha>0$, there exists $A_3^\partial=A_3^\partial(d, R, \alpha)\in (0, 1)$ such that the following statement holds.

    Given $t_0\in\R$, consider a Borel set $E\subseteq \partial C\cap (B_R(0_{\R^d})\times[t_0, t_0+R])$ with $\Haus^d(E)=1$.
    If $\varphi(t_0)\ge \alpha$ and $\abs{\varphi'(t_0)}\le A_3^\partial$, then
    \begin{equation*}
        \Per_{\partial C}(E) > I_{\partial\Sigma\times\R}(1).
    \end{equation*}
\end{lemma}
\begin{proof}
    It is proven repeating verbatim the proof of \cref{lem:crucial}, up to replacing $d$ with $d-1$, $\abs{\emptyparam}$ with $\Haus^{d}$, $C$ with $\partial C$, $\Sigma$ with $\partial\Sigma$.
\end{proof}

We are ready to choose the function $\varphi$ and show that $C$ and $\partial C$ may not have isoperimetric sets for small volumes.

\begin{proposition}\label{prop:main}
    Given $d\ge 2$, let $\Sigma\subseteq\R^d$ be a subset satisfying \cref{eq:sigma-properties}.

    There exists a function $\varphi:\R\to(0,\infty)$ satisfying \cref{eq:phi-properties} such that the convex set $C\subseteq\R^{d+1}$ defined in \cref{eq:def-C} does not have isoperimetric sets with volume equal to $1$.
\end{proposition}
\begin{proof}
    \noindent\textbf{Construction of $\varphi$.}
    We apply \cref{lem03} to construct $\varphi$ (so, in particular, $\varphi(0)=1$). As function $h$ in the lemma we take
    \begin{equation*}
        h(z) \defeq \frac12\sup_{0<z'\le z}A_3(d, A_1+A_2, z'),
    \end{equation*}
    where $A_3$ is defined in \cref{lem:crucial}, $A_2$ is defined in \cref{lem:Boundedness}, while $A_1 \defeq \sup_{x\in\Sigma} \frac{\abs{x}}{x\cdot e_1}$. Notice that $h$ is nondecreasing and strictly positive for $z>0$.
    We have defined $\varphi$ on $[0,\infty)$, we extend its domain to $\R$ by
    \begin{equation*}
        \varphi(t) \defeq \varphi(0) + t\varphi'(0) \quad\text{for $t<0$.}
    \end{equation*}
    The function $\varphi$ we have constructed satisfies all the properties in \cref{eq:phi-properties}.

    Assume by contradiction that $E\subseteq C$ with $\abs{E} = 1$ is an isoperimetric set, i.e., that $\Per_C(E) = I_{C}(\abs{E})$.  

    \vspace{2ex}
    \noindent\textbf{Space localization of the isoperimetric set $E$.}
    As a first step, we show that we can assume that $A_1 \ge \inf_{(x, t)\in E}\abs{x}$.
    Since $\varphi(0)=1$, in the region $\{(x, t): \abs{x}>A_1\}$ the set $C$ coincides with the cone $C'\subseteq\R^{d+1}$ defined in \cref{eq:outside-affine-cone}. Let $p'\in\R^{d+1}$ be the tip of such cone.
    If $\inf_{(x, t)\in E}\abs{x} > A_1$, we replace $E$ with a \emph{rescaled enlargement} of itself as follows. For $r>0$, denote with $E_r$ the set $\{z\in C: \mathrm{d}_{\overline E}(z)\le r\}$.
    We can find $r_0>0$ and $0<\lambda_0<1$ such that the set $E'$ defined as
    \begin{equation*}
        E' \defeq 
            p' + \lambda_0(E_{r_0}-p')
    \end{equation*}
    has volume $1$ and $\inf_{(x, t)\in E'}\abs{x}=A_1$. Thanks to \cref{eq:enlargement-ratio}, we have
    \begin{equation*}
        \frac{\Per_C(E')^{\frac{d+1}{d}}}{\abs{E'}}
        =
        \frac{\Per_C(E_{r_0})^{\frac{d+1}{d}}}{\abs{E_{r_0}}}
        \le
        \frac{\Per_C(E)^{\frac{d+1}{d}}}{\abs{E}},
    \end{equation*}
    so also $E'$ is an isoperimetric set. By replacing $E$ with $E'$, we gain the additional assumption $A_1\ge \inf_{(x,t)\in E}\abs{x}$.

    Applying \cref{eq:diam-C}, we know that $\diam(E)\le A_2$. In particular, for some $t_0\in\R$, it holds that $E\subseteq \{(x, t): t \in [t_0, t_0+A_2]\}$. Since we have reduced to the case $\inf_{(x, t)\in E}\abs{x} \le A_1$, the bound on the diameter implies that $\abs{x}\le A_2 + A_1$ for all $(x, t)\in E$. Hence, we have $E\subseteq C\cap (B_{A_1+A_2}(0_{\R^d})\times[t_0, t_0+A_2])$.

    \vspace{2ex}
    \noindent\textbf{Contradiction via \cref{lem:crucial}.}
    Both if $t_0<0$ or $t_0\ge 0$, we have that $\abs{\varphi'(t_0)}\le h(\varphi(t_0)) \le A_3(d, A_1+A_2, \alpha)$ for some $0<\alpha\le \varphi(t_0)$.
    Therefore we can apply \cref{lem:crucial} (with $R=A_1+A_2$) and deduce $\Per_C(E)>I_{\Sigma\times\R}(1)$. This is a contradiction because $\Sigma\times\R$ is a limit at infinity of $C$ and therefore $I_C(1)\le I_{\Sigma\times\R}(1)$.
\end{proof}

And we state the equivalent result on the boundary $\partial C$.

\begin{proposition}\label{prop:main2}
    Given $d\ge 3$, let $\Sigma\subseteq\R^d$ be a subset satisfying \cref{eq:sigma-properties}.

    There exists a function $\varphi:\R\to(0,\infty)$ satisfying \cref{eq:phi-properties} such that the hypersurface $\partial C\subseteq\R^{d+1}$ defined in \cref{eq:def-C} does not have isoperimetric sets with volume equal to $1$.
\end{proposition}
\begin{proof}
    It is proven repeating verbatim the proof of \cref{prop:main}, up to replacing $d$ with $d-1$, $\abs{\emptyparam}$ with $\Haus^{d}$, $C$ with $\partial C$, $\Sigma$ with $\partial\Sigma$. Moreover, one shall use \cref{lem:crucial-boundary} instead of \cref{lem:crucial} and \cref{eq:diam-C-boundary} instead of \cref{eq:diam-C}.
\end{proof}

\section{Smoothing the boundary of the construction}\label{sec:approximation}
We prove two simple lemmas (one for a convex set and one for its boundary) which show that the nonexistence of isoperimetric sets is stable under an appropriate approximation procedure. The proofs use crucially the results of \cite{AzagraStolyarov}.
\begin{proposition}\label{prop:smoothing}
    Given $d\ge 2$, let $\mathscr C\subseteq\R^d$ be a closed convex set such that $\partial \mathscr C$ does not contain any line (or, equivalently, $\mathscr C$ does not split a line) and its asymptotic cone $\mathscr C_\infty$ is nondegenerate. Assume that, for some $v>0$, any Borel set $E\subseteq \mathscr C$ with $\abs{E}=v$ satisfies $\Per_{\mathscr C}(E)>I_{\mathscr C}(v)$ (i.e., there is no isoperimetric set with volume $v$).

    There exists a strictly convex set $\tilde{\mathscr C}$ with smooth boundary so that\footnote{Here $\dist_H$ denotes the Hausdorff distance defined by using the Euclidean distance.} $\dist_H(\mathscr C\setminus B_R,\tilde{\mathscr C}\setminus B_R)\to 0$ as $R\to\infty$ and any Borel set $E\subseteq \tilde{\mathscr C}$ with $\abs{E}=v$ satisfies $\Per_{\tilde{\mathscr C}}(E)>I_{\tilde{\mathscr C}}(v)=I_{\mathscr C}(v)$.
\end{proposition}
\begin{proof}
    For any compact set $K\subseteq\R^d$, a simple compactness argument tells us that there exists $\eps_1=\eps_1(K, \mathscr C, v)>0$ such that if $E\subseteq \mathscr C\cap K$ has volume $\abs{E}=v$ then $\Per_{\mathscr C}(E)>I_{\mathscr C}(v) + \eps_1$. Therefore, with another compactness argument (which uses the lower semicontinuity of the perimeter along a sequence of Alexandrov spaces converging in the Gromov-Hausdorff sense~\cite[Proposition 3.6]{AmbrosioBrueSemola19}), we can find $\eps_2=\eps_2(K, {\mathscr C}, v)$ such that if $\tilde {\mathscr C}$ is a convex set with $\dist_H({\mathscr C}\cap K, \tilde {\mathscr C}\cap K)<\eps_2$ then any set $E\subseteq \tilde {\mathscr C}\cap K$ with $\abs{E}=v$ satisfies $\Per_{\tilde {\mathscr C}}(E)>I_{\mathscr C}(v)$.

    Let ${\mathscr C}_\infty$ be the asymptotic cone of $\mathscr C$. Thanks to \cref{lem:Boundedness}, we know that there exists a constant $L=L({\mathscr C}_\infty,v)>0$ such that if $\tilde {\mathscr C}$ is a convex set with asymptotic cone ${\mathscr C}_\infty$, then any isoperimetric set in $\tilde {\mathscr C}$ with volume $\le v$ has diameter below $L$.
    
    Let $(K_n)_{n\in\N}$ be a sequence of compact sets such that $\cup_{n\in\N}K_n=\R^d$ is a \emph{locally finite} covering, and for any set $E\subseteq\R^d$ with diameter below $L$ there is $n_0=n_0(E)$ such that $E\subseteq K_{n_0}$. For every $n\in\mathbb N$, let us define
    \begin{equation}
    U_{n}:=\left(\mathbb R^d\setminus K_n\right)\cup\{x\in\mathbb R^d:\mathrm{d}_{\partial {\mathscr C}\cap K_n}(x)< \min\{\varepsilon_2(K_n,{\mathscr C},v),n^{-1}\}\}.
    \end{equation}
    Moreover, let us define
    \begin{equation}
        U\defeq \bigcap_{n\in\mathbb N} U_{n}.
    \end{equation}
    Applying \cite[Theorem 2]{AzagraStolyarov} to the convex set ${\mathscr C}$, using $U$ as  open set, we get that there exists a closed strictly convex set $\tilde {\mathscr C}$ with nonempty interior and smooth boundary such that 
    \begin{itemize}
        \item $\dist_H({\mathscr C}\setminus B_R,\tilde {\mathscr C}\setminus B_R)\to 0$ as $R\to +\infty$.
        \item $\dist_H({\mathscr C}\cap K_n, \tilde {\mathscr C}\cap K_n)<\eps_2(K_n, {\mathscr C}, v)$.
    \end{itemize}

    Notice that the pointed limits at infinity of ${\mathscr C}$ coincide with the pointed limits at infinity of $\tilde {\mathscr C}$ and the asymptotic cone of ${\mathscr C}$ coincides with the asymptotic cone of $\tilde {\mathscr C}$.
    
    Since $I_{\mathscr C}(v)$ is not achieved on ${\mathscr C}$, there is a pointed limit at infinity of ${\mathscr C}$ where there is an isoperimetric set with volume $v$ and perimeter $I_{\mathscr C}(v)$, see \cref{lem:LimitiVannoAInfinito}.
   Since also $\tilde {\mathscr C}$ has such pointed limit at infinity, we deduce, by \cref{prop:profile-properties}, that $I_{\tilde {\mathscr C}}(v)\le I_{\mathscr C}(v)$.

    Assume by contradiction that there is a set $E\subseteq \tilde {\mathscr C}$ with $\abs{E}=v$ and $\Per_{\tilde {\mathscr C}}(E)=I_{\tilde {\mathscr C}}(v)$. Then the diameter of $E$ is bounded above by $L$ and therefore, by construction, there is a compact set $K_{n_0}$ such that $E\subseteq K_{n_0}$. Then, by construction of $\tilde {\mathscr C}$ and by definition of $\eps_2(K_{n_0}, {\mathscr C}, v)$, we deduce that $I_{\tilde {\mathscr C}}(v)=\Per_{\tilde {\mathscr C}}(E)>I_{\mathscr C}(v)$ which is a contradiction.

    To conclude, observe that since $I_{\tilde {\mathscr C}}(v)$ is not achieved on $\tilde {\mathscr C}$, it must be achieved on a pointed limit at infinity of $\tilde {\mathscr C}$. And arguing as before we deduce $I_{\tilde {\mathscr C}}(v)\ge I_{\mathscr C}(v)$.
\end{proof}

\begin{proposition}\label{prop:smoothing-boundary}
    Given $d\ge 3$, let ${\mathscr C}\subseteq\R^d$ be a convex set such that $\partial {\mathscr C}$ does not contain any line (or, equivalently, ${\mathscr C}$ does not split a line) and its asymptotic cone ${\mathscr C}_\infty$ is nondegenerate. Assume that, for some $v>0$, any Borel set $E\subseteq \partial {\mathscr C}$ with $\Haus^{d-1}(E)=v$ satisfies $\Per_{\partial \mathscr C}(E)>I_{\partial {\mathscr C}}(v)$ (i.e., there is no isoperimetric set with volume $v$).

    There exists a strictly convex set $\tilde {\mathscr C}$ with smooth boundary so that $\dist_H(\partial {\mathscr C}\setminus B_R, \partial \tilde {\mathscr C}\setminus B_R)\to 0$ as $R\to\infty$ and any Borel set $E\subseteq  \partial \tilde{\mathscr C}$ with $\Haus^{d-1}(E)=v$ satisfies $\Per_{\partial\tilde{\mathscr C}}(E)>I_{\partial \tilde {\mathscr C}}(v)=I_{\partial {\mathscr C}}(v)$.
\end{proposition}
\begin{proof}
    It is proven repeating verbatim the proof of \cref{prop:smoothing}, up to replacing $d$ with $d-1$, $\abs{\emptyparam}$ with $\Haus^{d-1}$, ${\mathscr C}$ with $\partial {\mathscr C}$, $\tilde {\mathscr C}$ with $\partial \tilde {\mathscr C}$ appropriately.
\end{proof}

\section{Proof of the main results}\label{sec:proofs}

\begin{proof}[Proof of \cref{thm:Main}]
    For $d\geq 3$, let $C\subseteq \mathbb R^d$ be the convex set constructed  in \cref{prop:main}. Thus $C$ does not have isoperimetric sets of volume $1$. By \cref{prop:smoothing} we can construct a strictly convex set $\tilde C\subseteq\mathbb R^{d}$ with smooth boundary such that $\tilde C$ does not have isoperimetric sets of volume $1$. Notice that, thanks to the construction of $\tilde C$ in \cref{prop:smoothing}, $\tilde C$ has the same asymptotic cone and the same limits at infinity of $C$. Since $C$ has nondegenerate asymptotic cone and all its limits at infinity are metric cones, see \cref{PropertiesC}, we can apply \cref{lem:Threshold} to $\tilde C$. Since $\tilde C$ does not have isoperimetric sets of volume $1$, we get that the threshold $v_0$ coming from \cref{lem:Threshold} verifies $v_0>0$. After properly scaling in such a way $v_0=1$, the theorem is shown with $C=\tilde C$.
\end{proof}

\begin{proof}[Proof of \cref{thm:Main-boundary}]
    For $d\geq 3$, by \cref{prop:main2} there exists a convex set $C\subseteq\mathbb R^{d+1}$ such that $\partial C\subseteq\mathbb R^{d+1}$ does not have isoperimetric sets of volume $1$. By \cref{prop:smoothing-boundary} we can construct a strictly convex set $\tilde C\subseteq\mathbb R^{d+1}$ with smooth boundary $\partial \tilde C\subseteq\mathbb R^{d+1}$ such that $\partial \tilde C$ does not have isoperimetric sets of volume $1$. Notice that, thanks to the properties of $\tilde C$ stated in \cref{prop:smoothing-boundary}, $\partial \tilde C$ has the same asymptotic cone and the same limits at infinity of $\partial C$. Observe that, by construction of $C$ (see \cref{subsec:presentation-construction}), all the limits at infinity of $\partial C$ are cones, and $\partial C$ has nondegenerate asymptotic cone. 
    We are in position to apply \cref{lem:Threshold} to $\partial \tilde C$, and then obtain a threshold $v_0>0$. After scaling in such a way that $v_0=1$, we get the sought conclusion by taking $M=\partial\tilde C$.
\end{proof}

\appendix
\section{A simple technical lemma} \label{app:boring}
The following standard technical lemma was needed in the proof of \cref{prop:main}.
\begin{lemma}\label{lem03}
    Let $h:(0,1]\to(0,\infty)$ be a (weakly) increasing function.
    There exists a continuously differentiable convex (strictly) decreasing function $\varphi:[0,\infty)\to (0, 1]$ such that $\varphi(0)=1$, $\varphi(+\infty)=0$, and $\abs{\varphi'(t)} \le h(\varphi(t))$ for all $t\ge 0$.
\end{lemma}
\begin{proof}
    As a first step, let $\tilde h:[0, 1]\to[0,\infty)$ be a strictly increasing function that is Lipschitz continuous and satisfies $\tilde h\le h$ on $(0, 1]$ and $\tilde h(0)=0$. For example, it is enough to choose $\tilde h(t):=\int_0^t h(s)\mathrm{d}s$.

    Consider the Cauchy problem
    \begin{align}
        &\varphi(0) = 1 \notag, \\
        &\varphi'(t) = -\tilde h(\varphi(t)). \label{eq:49}
    \end{align}
    Observe that the constant function $0$ is a solution of \cref{eq:49}, while the constant function $1$ is a super-solution of \cref{eq:49}.
    Hence there exists a unique solution $\varphi:[0,\infty)\to(0,1]$ of the Cauchy problem (local existence follows from \cite[Chapter 1, Theorem 2.3]{CoddingtonLevinson}, and global existence from the presence of the two mentioned barriers~\cite[Chapter 1, Theorem 4.1]{CoddingtonLevinson}).

    Since the function $\varphi$ satisfies \cref{eq:49}, we get that $\varphi'<0$ and so $\varphi$ is decreasing. Thus the quantity $-\tilde h(\varphi(t))$ is increasing and therefore $\varphi'$ is increasing, which is equivalent to the convexity of $\varphi$. Also $\abs{\varphi'(t)} \le h(\varphi(t))$ follows from \cref{eq:49}.

    It remains only to see that $\varphi(+\infty)=0$. Since $\varphi$ is decreasing, it has a (nonnegative) limit; assume by contradiction that $\varphi(+\infty)=\ell>0$. Then,
    \begin{equation*}
        \varphi'(t) = -\tilde h(\varphi(t)) \le -\tilde h(\ell),
    \end{equation*}
    for all $t\ge 0$ and therefore
    \begin{equation*}
        \varphi(t) \le \varphi(0) - \tilde h(\ell)t.
    \end{equation*}
    This produces the desired contradiction as it implies that $\varphi(t)<0$ for large values of $t$ (since $\tilde h(\ell)>\tilde h(0) = 0$).
\end{proof}

\printbibliography

\end{document}